\documentclass[aps,reprint,twocolumn,showkeys,superscriptaddress,longbibliography]{revtex4-2}
\usepackage{bm}
\usepackage{lipsum}
\usepackage{graphicx}
\usepackage{algorithm}
\usepackage{algpseudocode}

\usepackage{amsmath}
\usepackage{hyperref} 
\hypersetup{
colorlinks=true,
citecolor=blue, 
urlcolor=blue,
linkcolor=blue
}
\usepackage[symbol]{footmisc}

\begin{document}
\date{\today}
\title{Community detection by simulated bifurcation}

\author{Wei Li}
\email{qstatex@gmail.com}
\noaffiliation

\author{Yi-Lun Du}
\email{yilun.du@iat.cn}
\affiliation{Shandong Institute of Advanced Technology,  Jinan, China}

\author{Nan Su}
\email{banban.su@gmail.com}
\noaffiliation

\author{Konrad Tywoniuk}
\affiliation{Department of Physics and Technology, University of Bergen, Bergen, Norway}

\author{Kyle Godbey}
\affiliation{Facility for Rare Isotope Beams, Michigan State University, Michigan, USA}

\author{Horst St\"ocker}
\affiliation{Frankfurt Institute for Advanced Studies, Frankfurt am Main, Germany}
\affiliation{Institut f\"ur Theoretische Physik, Goethe Universit\"at, Frankfurt am Main, Germany}
\affiliation{GSI Helmholtzzentrum f\"ur Schwerionenforschung GmbH, Darmstadt, Germany}

\begin{abstract}
Community detection, also known as graph partitioning, is a well-known NP-hard combinatorial optimization problem with applications in diverse fields such as complex network theory, transportation, and smart power grids. The problem's solution space grows drastically with the number of vertices and subgroups, making efficient algorithms crucial. In recent years, quantum computing has emerged as a promising approach to tackling NP-hard problems. This study explores the use of a quantum-inspired algorithm, Simulated Bifurcation (SB), for community detection. Modularity is employed as both the objective function and a metric to evaluate the solutions. The community detection problem is formulated as a Quadratic Unconstrained Binary Optimization (QUBO) problem, enabling seamless integration with the SB algorithm. Experimental results demonstrate that SB effectively identifies community structures in benchmark networks such as Zachary's Karate Club and the IEEE 33-bus system. Remarkably, SB achieved the highest modularity, matching the performance of Fujitsu's Digital Annealer, while surpassing results obtained from two quantum machines, D-Wave and IBM. These findings highlight the potential of Simulated Bifurcation as a powerful tool for solving community detection problems.

\end{abstract}

\maketitle

\section{Introduction}

Networks provide a powerful framework for modeling the structure and dynamics of complex systems by representing them as nodes connected by edges~\cite{boccaletti2006complex,newman2018networks}. This framework finds broad application across diverse fields, including social networks~\cite{mislove2011understanding,du2007community,kojaku2021effectiveness}, transportation systems~\cite{colizza2006role,barthelemy2011spatial,lin2013complex}, finance~\cite{caldarelli2004emergence,bardoscia2021physics,barucca2020network}, and neuroscience~\cite{bassett2006small,sporns2014contributions,bassett2017network}. The modern science of networks seeks to unravel both the structural and functional aspects of these systems, where nodes represent fundamental units and edges denote their interactions~\cite{fortunato202220}.

One of the most critical features of many networks is their community structure—clusters of nodes that are more densely connected internally than with the rest of the network. Community detection, also referred to as graph partitioning, is a central challenge in network science, aiming to divide the network into groups while minimizing inter-group connections~\cite{girvan2002community}. This problem is classified as NP-hard, making it increasingly difficult to solve as network size grows~\cite{lucas2014ising}.

Among the various techniques for community detection, modularity maximization is one of the most widely used. Modularity is a scalar function that measures the quality of a network division by comparing the density of links within communities against a random baseline~\cite{newman2004finding}. A higher modularity score indicates a more pronounced community structure. However, maximizing modularity becomes computationally infeasible for large networks, necessitating the use of approximate methods.

Several heuristic approaches have been developed to address modularity maximization, including greedy algorithms~\cite{newman2004fast,clauset2004finding}, simulated annealing~\cite{guimera2004modularity,medus2005detection}, extremal optimization~\cite{duch2005community}, genetic algorithms~\cite{li2010genetic}, and the widely used Louvain method~\cite{blondel2008fast}. Despite these advances, no single algorithm performs optimally across all network types due to variations in network structures and purposes~\cite{rostami2023community}. The vast search space further complicates solving community detection problems with conventional methods~\cite{mohseni2022ising}.

To overcome these challenges, Ising computing has emerged as a promising approach for solving hard combinatorial optimization problems, including community detection~\cite{mohseni2022ising}. The Ising model is particularly notable because any NP-hard problem can be reformulated as an Ising model with polynomial overhead. Scalable Ising computing frameworks could revolutionize optimization in various domains. Among the leading implementations of Ising computing are hardware solutions like D-Wave’s Quantum Annealer (QA) and Fujitsu’s Digital Annealer (DA), both of which have demonstrated high efficiency in tackling combinatorial optimization tasks. Recent studies have applied these systems to community detection with encouraging results~\cite{negre2020detecting,fernandez2021community,gemeinhardt2021quantum,Wierzbinski2022brain,kao2023solving}.

In addition to hardware solutions, Ising computing can also be implemented algorithmically on conventional hardware. One such approach is the Simulated Bifurcation (SB) algorithm, a GPU-accelerated quantum-inspired method designed for solving large-scale combinatorial problems efficiently~\cite{goto2016bifurcation}. This study explores the application of SB to community detection, benchmarking its performance on two networks: Zachary’s Karate Club, a well-known social network, and the IEEE 33-Bus system, a benchmark in power distribution networks. The results are compared with those of state-of-the-art quantum and quantum-inspired computers, including IBM’s quantum systems, D-Wave’s Quantum Annealer and Fujitsu’s Digital Annealer, as well as mathematical optimization software Gurobi. Surprisingly, our results outperform the IBM and D-Wave quantum computers and Gurobi, and achieve the same accuracy as the Fujitsu quantum-inspired computer.

\section{Modularity-Based Community Detection\label{sec2}}

This section presents the theoretical framework for community detection in complex networks, describes the modularity maximization problem, and explains how it is transformed into a quadratic unconstrained binary optimization (QUBO) formulation suitable for quantum or quantum-inspired computation.

\subsection{Modularity Formulation}

Community detection involves identifying cohesive subgraphs within a network \( G(V, E) \), where \( V \) and \( E \) represent the sets of vertices and edges, respectively. The objective is to uncover large-scale patterns that emerge from interactions among individual nodes. A key metric used to evaluate the quality of community structures is \textit{modularity} \( Q_e \)~\cite{newman2004finding}. For a weighted graph, the modularity score is defined as:
\begin{equation}
Q_e = \frac{1}{2m} \sum_{i,j} \left( A_{ij} - \gamma \frac{k_i k_j}{2m} \right) \delta(c_i, c_j),
\label{eq: Qe}
\end{equation}
where:
\begin{itemize}
    \item \( A_{ij} \): the element in row \( i \) and column \( j \) of the adjacency matrix of \( G \),
    \item \( m = \frac{1}{2} \sum_{i,j} A_{ij} \): the total weighted number of edges,
    \item \( k_i = \sum_{j} A_{ij} \): the degree of node \( i \), i.e., sum of weighted edges connected to node \( i \),
    \item \( \delta(c_i, c_j) \): a Kronecker delta function that equals 1 if nodes \( i \) and \( j \) belong to the \textit{same community}, and 0 otherwise,
    \item \( \gamma \): a resolution parameter, set to 1 in this study.
\end{itemize}

The first term in \( Q_e \) measures the proportion of edge weights within communities, while the second term estimates the expected fraction of edge weights under a null model where connections are random. Higher values of modularity indicate stronger intra-community connectivity compared to random configurations.

The total number of possible community partitions for a graph with \( n \) nodes, where $K\leq n$, is given by the \textit{Bell number}:
\begin{equation}
B_n = \frac{1}{e} \sum_{m=0}^\infty \frac{m^n}{m!},
\end{equation}
whose asymptotic upper bound is~\cite{Daniel2010Bell}:
\begin{equation}
B_n \leq \left( \frac{0.792n}{\ln(n+1)} \right)^n.
\label{eq:B_max}
\end{equation}

Finding the optimal partition into \( K_{\text{opt}} \) communities involves maximizing the modularity \( Q_e \) as a function of \( K \).

\subsection{QUBO Formulation}

To leverage quantum or quantum-inspired computing for community detection, the modularity function must be mapped to a mathematical framework of Ising spin glass model.

\subsubsection{Ising Model}

The Ising spin glass model minimizes the spin system energy given by:
\begin{equation}
E(S) = - \sum_{i=1}^N \sum_{j=1}^N J_{i,j} s_i s_j - \sum_{i=1}^N h_i s_i,
\end{equation}
where:
\begin{itemize}
    \item \( s_i \in \{-1, 1\} \): the binary spin state of the \( i \)-th spin,
    \item \( J_{i,j} \): the coupling coefficient between spins \( i \) and \( j \), satisfying \( J_{i,j} = J_{j,i} \) and \( J_{i,i} = 0 \),
    \item \( h_i \): the external magnetic field acting on spin \( i \),
    \item \( N \): the total number of spins.
\end{itemize}

\subsubsection{Mapping to QUBO}

The Ising model can be transformed into a QUBO problem using the substitution \( s_i = 2x_i - 1 \), where \( x_i \in \{0, 1\} \). This reformulation yields:
\begin{equation}
H(x) = \bm{x}^T \hat{Q} \bm{x},
\end{equation}
where:
\begin{itemize}
    \item \( \bm{x} \): the binary vector of variables \( x_i \),
    \item \( \hat{Q} \): a symmetric matrix whose elements derive from the Ising model parameters.
\end{itemize}

\subsubsection{Community Detection as a QUBO Problem}

In the context of community detection, each node \( i \) is assigned a binary vector \( \bm{\mathrm{x}}_i = (x_{i0}, x_{i1}, \ldots, x_{i(K-1)}) \), where \( K \) is the number of communities. If node \( i \) belongs to community \( k \), then \( x_{ik} = 1 \) and all other entries in \( \bm{\mathrm{x}}_i \) are 0. The modularity function to be maximized is reformulated as minimizing $M$:
\begin{equation}
M = - \bm{x}^T \hat{Q}_e \bm{x},
\label{eq5}
\end{equation}
where \( \bm{x} = (x_{00}, x_{01}, \ldots, x_{n-1, K-1})^T \) is a flattened binary vector of size $nK$, and \( \hat{Q}_e \) encodes the modularity matrix according to Eq.~(\ref{eq: Qe}).

\subsubsection{Constraints}

Two constraints must be enforced to ensure valid community assignments:
\begin{enumerate}
    \item \textbf{Each node belongs to exactly one community}:
    \begin{equation}
    \sum_{k=0}^{K-1} x_{ik} = 1, \quad \text{for } i = 0, 1, \ldots, n-1.
    \end{equation}
    \item \textbf{Each community contains at least one node}:
    \begin{equation}
    \sum_{i=0}^{n-1} x_{ik} \geq 1, \quad \text{for } k = 0, 1, \ldots, K-1.
    \end{equation}
\end{enumerate}

These constraints are incorporated into the Hamiltonian using penalty terms. The complete Hamiltonian is given by:
\begin{widetext}
\begin{equation}
H \equiv  - {\bm{x}'}^T \hat{Q}_e' \bm{x}' = - \bm{x}^T \hat{Q}_e \bm{x} + \alpha \sum_{i=0}^{n-1} \left( \sum_{k=0}^{K-1} x_{ik} - 1 \right)^2 
+ \beta \sum_{k=0}^{K-1} \left( \sum_{i=0}^{n-1} x_{ik} - \sum_{d=1}^{d_{max}} 2^{d-1} x_{dk} - 1 \right)^2,
\label{eq: H}
\end{equation}
\end{widetext}
where \( \alpha \) and \( \beta \) are the penalty multipliers for enforcing the constraints. To satisfy the second constraint, binary slack variables 
$x_{dk}$ are introduced, forming a new binary vector $\bm{x}'$  by augmenting $\bm{x}$.  $\sum_{i=0}^{n-1} x_{ik}$ is enforced to equal the term $\sum_{d=1}^{d_{max}} 2^{d-1} x_{dk}+1$, which can represent any integer number from $1$ to $2^{d_{max}}$ and provides the flexibility to accommodate different community configurations. The value of $d_{max}$ is predetermined based on the number of nodes and the communities in the network. Both the slack variables $x_{dk}$ and the primary decision variables $x_{ik}$ are optimized simultaneously by the Simulated Bifurcation (SB) algorithm to ensure high modularity while satisfying the constraints.

The resulting $\hat{Q}_e'$ matrix in Eq.~(\ref{eq: H}), serves as the input for the SB algorithm. To ensure optimal performance, we carefully configure the penalty constants in Eq.~(\ref{eq: H}), maximizing modularity while enforcing the problem's constraints. This reformulation bridges modularity-based community detection with the QUBO framework, enabling efficient exploration of large solution spaces using quantum or quantum-inspired algorithms.

\section{Methodology}
To solve the problem of modularity-based community detection, this study adopts the Simulated Bifurcation (SB) algorithm~\cite{goto2016bifurcation}, which is well suited to operate on GPUs. The SB algorithm relies on an adiabatic and ergodic search mechanism for the Ising problem. Furthermore, two enhanced variants of SB, namely the ballistic Simulated Bifurcation (bSB) and discrete Simulated Bifurcation (dSB) algorithms~\cite{goto2021high}, have been developed. These variants offer improvements in both computational speed and solution accuracy. As presented in~\cite{kanao2022simulated}, both bSB and dSB are based on the equations of motion of following Hamiltonian $H_{\textrm{SB}}$,
\begin{flalign}
        \dot{x}_i =& \frac{\partial H_{\textrm{SB}}}{\partial y_i} = a_0 y_i, \\
        \dot{y}_i =& - \frac{\partial H_{\textrm{SB}}}{\partial x_i} = - \left[ a_0 - a(t) \right] x_i + c_0 f_i,  \label{eq:dSB} \\
        H_{\textrm{SB}} =& \frac{a_0}{2} \sum_{i=1}^{N} y_i^2 + V_{\textrm{SB}},
\end{flalign}
where $x_i$ and $y_i$ represent respectively the \textit{positions} and \textit{momenta} corresponding to the spin state $s_i$ in this section, while the dots denotes time derivatives. The control parameter $a(t)$ is a time-dependent variable increasing monotonically from zero at the initial time to a positive constant $a_0$ at the final time, and $c_0$ is another positive constant. The potential $V_{\textrm{SB}}$ and interaction forces $f_i$ differ between the bSB and dSB formulations, defined as follows:
\begin{equation}
V_{\textrm{SB}} = 
\begin{cases} 
\begin{aligned}
    & \frac{a_0 - a(t)}{2} \sum_{i=1}^{N} x_i^2 - \frac{c_0}{2} \sum_{i=1}^{N} \sum_{j=1}^{N} J_{i,j} x_i x_j, \\
    & \qquad\text{when } |x_i| \leq 1 \text{ for all } x_i \text{ (for bSB)}, \\
    & \frac{a_0 - a(t)}{2} \sum_{i=1}^{N} x_i^2 - \frac{c_0}{2} \sum_{i=1}^{N} \sum_{j=1}^{N} J_{i,j} x_i \text{sgn}(x_j), \\
    & \qquad\text{when } |x_i| \leq 1 \text{ for all } x_i \text{ (for dSB)}, \\
    & \infty, \qquad\text{otherwise},
\end{aligned}
\end{cases}
\end{equation}
\begin{equation}
 f_i = 
\begin{cases} 
\begin{aligned}
    & \sum_{j=1}^NJ_{i,j}x_j \quad \text{for bSB}, \\
    & \sum_{j=1}^NJ_{i,j} \mathrm{sgn}(x_j) \quad \text{for dSB},
\end{aligned}
\end{cases}
\end{equation}
where $ \mathrm{sgn}(x_j)$ denotes the sign of $x_j$. In this study, we focus on the dSB algorithm. By solving the equations of motion on GPUs using the \textit{symplectic Euler method}~\cite{leimkuhler2004simulating}, the algorithm ensures that a solution corresponding to at least a local minimum of the Ising problem is obtained at the final time. The spin states \( s_i \) are determined as \( s_i = \text{sgn}(x_i) \). The SB algorithm also allows for solving the QUBO problem directly with appropriate transformations. These resulting binary variables will represent the community assignments, effectively revealing the underlying community structure of the graph.

\section{Results and discussion}
This section illustrates the results of modularity-based community detection using simulated bifurcation algorithm on a single GPU, applied to two benchmark datasets.

\subsection{Karate Club}
\begin{figure}[tbh!]
\includegraphics[width=0.48\textwidth]{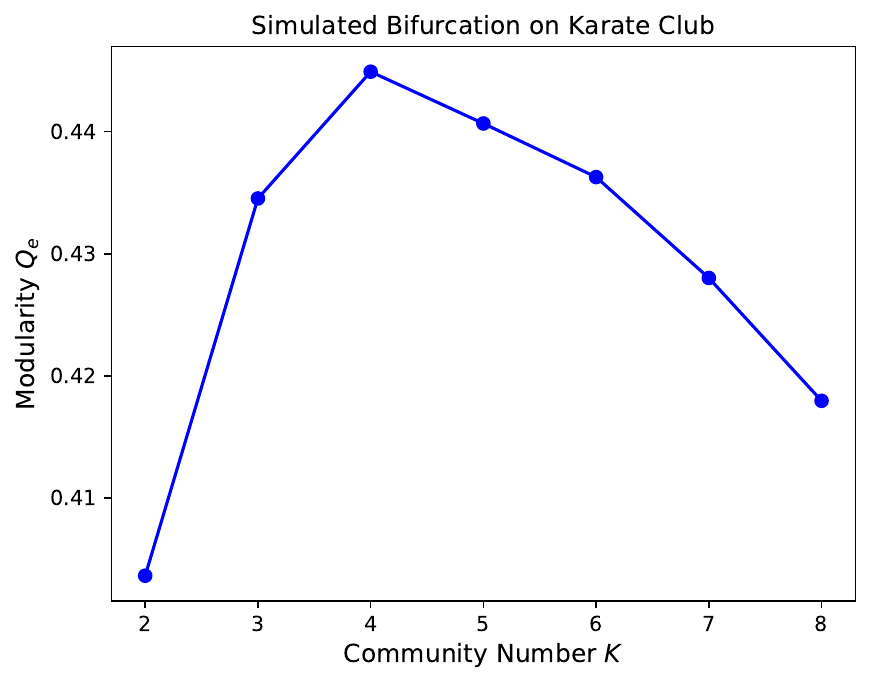}
\caption{Modularity ($Q_e$) versus Community Number ($K$) for Zachary's Karate Club using Simulated Bifurcation. The optimal result is achieved at $K_{\mathrm{opt}} = 4$ with $Q_e = 0.445$.}
\label{fig:Karate_modularity} 
\end{figure}

\begin{table}[tbh!]%[H] add [H] placement to break table across pages
       \begin{ruledtabular}
            \begin{tabular}{lccccc}
             Method &IBM~\cite{gemeinhardt2021quantum} &D-Wave~\cite{Wierzbinski2022brain} &Fujitsu~\cite{kao2023solving}& SB\\
             $Q_e$ &0.420 &0.444 &0.445  &0.445
            \end{tabular}
    \end{ruledtabular}
     \caption{Comparison of Best Community Detection Results on Zachary's Karate Club Across Different Hardware Platforms and Algorithms.\label{table1}}
\end{table}

Zachary's Karate Club network~\cite{zachary1977information} represents social interactions among 34 members of a university karate club. The network data was imported using the Python library NetworkX~\cite{hagberg2008exploring}. According to Eq.~(\ref{eq:B_max}), the upper bound of the total number of feasible solutions for this problem is approximately $7.89 \times 10^{29}$, highlighting the vastness of the solution space. 

\begin{figure}[tbh!]
    \includegraphics[width=0.5\textwidth]{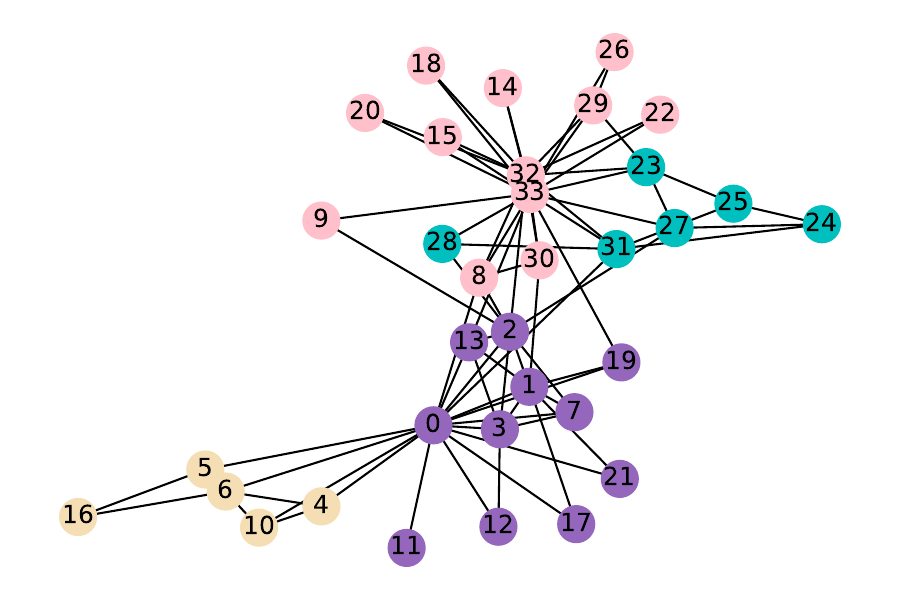}
    \caption{Partitioned Graph with the highest Modularity $Q_e=0.445$, with Community Number $K_{\mathrm{opt}}=4$ by Simulated Bifurcation for Zachary's Karate Club Network. \label{fig:Karate_best}}
\end{figure}

Fig.~\ref{fig:Karate_modularity} illustrates the modularity \( Q_e \) as a function of the number of communities \( K \), obtained using the SB algorithm. The SB algorithm identifies the optimal community structure at \( K_{\mathrm{opt}} = 4 \), achieving a modularity score of \( Q_e = 0.445 \). Tab.~\ref{table1} compares the best results from different hardware platforms and algorithms. Remarkably, our SB algorithm achieves a modularity of 0.445, outperforming IBM’s result of 0.419~\cite{gemeinhardt2021quantum} and D-Wave’s result of 0.444~\cite{Wierzbinski2022brain}. Moreover, it matches the accuracy of Fujitsu’s Digital Annealer~\cite{kao2023solving}. The partitioned graph corresponding to the highest modularity \( Q_e \) is shown in Fig.~\ref{fig:Karate_best}, which aligns perfectly with the result obtained by Fujitsu's method.

\subsection{Electrical Virtual Micriogrids}

With the rise of renewable energy, conventional centralized power distribution networks must adapt to diverse methods of electricity generation and transmission. Virtual microgrids have emerged as a promising solution to address this challenge, offering benefits such as reduced power loss during transmission, higher efficiency, and better compatibility with green energy sources~\cite{xu2018structural}. In this study, we apply the SB algorithm on the well-known IEEE 33-bus system benchmark~\cite{baran1989network}, imported from the Python package PandaPower~\cite{thurner2018pandapower} and converted into NetworkX graphs using the “create\_\,nxgraph” function. According to Eq.~(\ref{eq:B_max}), the upper bound of the total number of feasible solutions for this problem is approximately $5.09 \times 10^{28}$, which indicates the large size of the search space.

%%%%%%%%%%%%%%%%%%%%%%%%%%%%%%%%%%%%%%%%%%%%%%%%%%
\begin{figure}[tbh!]
\includegraphics[width=0.48\textwidth]{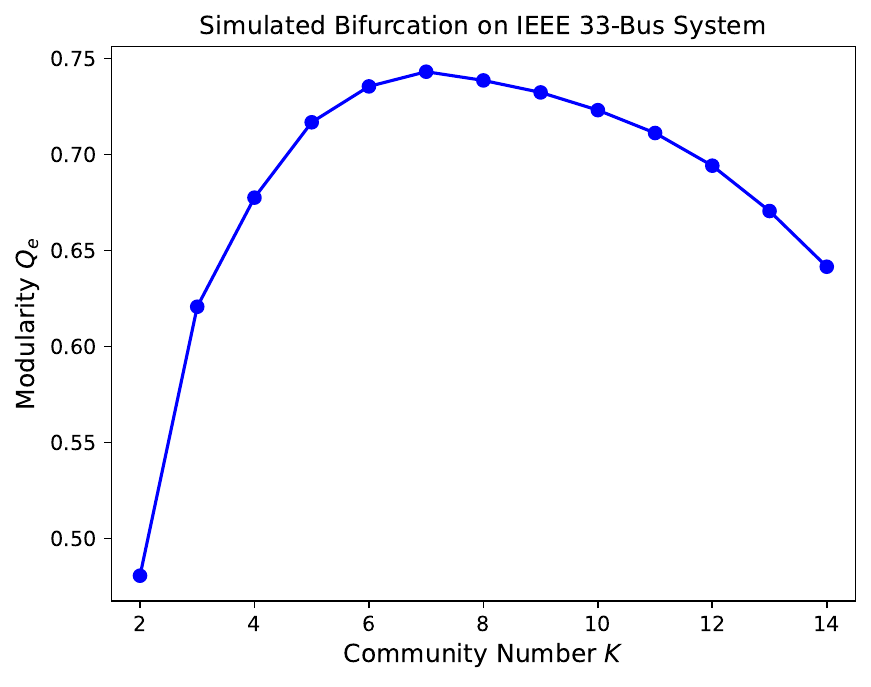}
\caption{Modularity ($Q_e$) versus Community Number ($K$) for the IEEE 33–Bus System using Simulated Bifurcation. The optimal result is achieved at $K_{\mathrm{opt}}=7$ with $Q_e=0.743$.}
\label{fig:IEEE33_modularity} 
\end{figure}
%%%%%%%%%%%%%%%%%%%%%%%%%%%%%%%%%%%%%%%%%%%%%%%%%%

For modularity-based graph partitioning applied to virtual microgrids in power distribution networks, we adopt the electrical modularity proposed by Kao et al.~\cite{kao2023solving}. In this framework, the weights of the edges are determined by the inverse of the absolute value of impedance~\cite{cotilla2013multi}. Specifically, given the resistance $r$ and reactance $x$ of each edge, the weights are calculated as $1 / | r + \dot{\imath} x |$. This approach is then used to compute the modularity of the network.

\begin{figure}[tbh!]
    \includegraphics[width=0.5\textwidth]{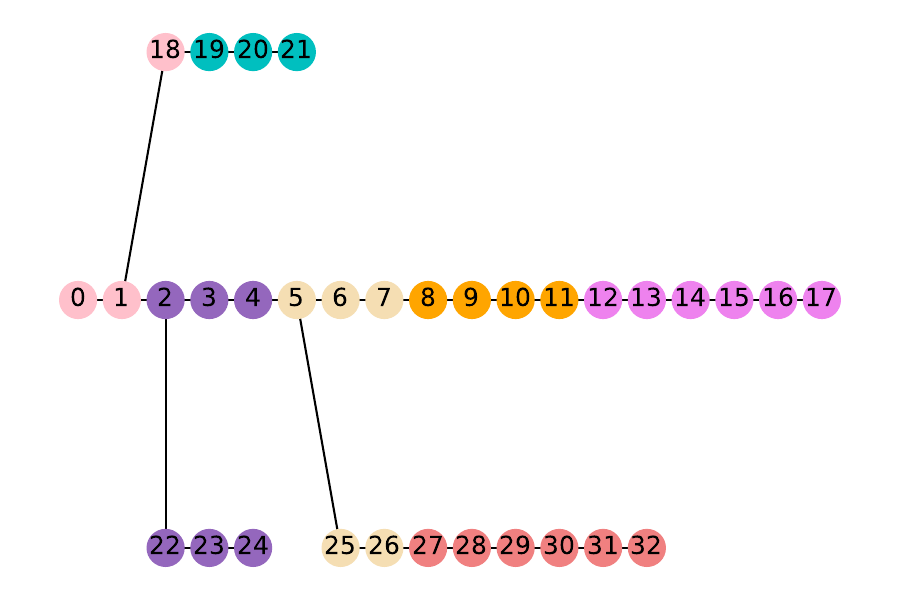}
    \caption{Partitioned Graph with the highest Modularity $Q_e=0.743$, with Community Number $K_{\mathrm{opt}}=7$ by Simulated Bifurcation for IEEE 33-Bus System. \label{fig:IEEE33_best}}
\end{figure}

The SB results of modularity versus community number for the IEEE 33-bus system are shown in Fig.~\ref{fig:IEEE33_modularity}, and the optimized results from different hardware platforms and algorithms are summarized in Tab.~\ref{table2}. Similar to the Karate Club case, our SB optimal result of 0.743 achieved at $K_\mathrm{opt}=7$ surpasses both the Gurobi result of 0.711~\cite{fernandez2021community} and the D-Wave result of 0.711~\cite{fernandez2021community}, and again achieves the same accuracy as Fujitsu’s result~\cite{kao2023solving}. The partitioned graph with the highest modularity \( Q_e = 0.743 \) is shown in Fig.~\ref{fig:IEEE33_best}, which corresponds to the same partition obtained by Fujitsu's method.

\begin{table}[tbh!]
      \begin{ruledtabular}
            \begin{tabular}{lcccc}
             Method &Gurobi~\cite{fernandez2021community} &D-Wave~\cite{fernandez2021community} &Fujitsu~\cite{kao2023solving} &SB\\
             $Q_e$ &0.711  &0.711  &0.743  &0.743
            \end{tabular}
    \end{ruledtabular}
    \caption{Comparison of Best Community Detection Results on IEEE 33-Bus System Across Different Hardware Platforms and Algorithms.\label{table2}}
\end{table}

\section{Conclusion}

In this manuscript, we studied the task of community detection using for the first time the technique of simulated bifurcation. The benchmarking problems we chose are computationally demanding for conventional techniques, and have been studied lately by quantum and quantum-inspired computers to show their effectiveness compared to conventional techniques. Our results demonstrate that the quantum-inspired SB algorithm, running on a single GPU, can outperform the IBM and D-Wave quantum computers, as well as a leading mathematical optimization software, Gurobi. Furthermore, SB matches the precision achieved by the quantum-inspired Fujitsu Digital Annealer. This is a highly nontrivial finding, especially considering the significantly reduced computational resources required by SB compared to the costly quantum computers. Additionally, GPU-accelerated SB offers greater flexibility in addressing computationally demanding problems, presenting an alternative to quantum hardware. As quantum computers remain constrained in size and capability due to technical challenges, quantum-inspired algorithms like SB hold timely significance for solving practical problems in science and industry.

\begin{acknowledgments}
We would like to acknowledge the support of the cluster Goethe-HLR of the Center for Scientific Computing. This work is supported by the Taishan Scholars Program under Grant No. tsqnz20221162, Shandong Excellent Young Scientists Fund Program (Overseas) under Grant No. 2023HWYQ-106 (Y. D.).
\end{acknowledgments}

\bibliography{main}

\end{document}